\documentclass[12TP,draft]{article}

\begin{document}

\begin{center}
\LARGE\noindent\textbf{A sufficient condition for  pre-Hamiltonian cycles in  bipartite digraphs }\\

\end{center}
\begin{center}
\noindent\textbf{Samvel Kh. Darbinyan and Iskandar A. Karapetyan}\\

Institute for Informatics and Automation Problems, Armenian National Academy of Sciences

E-mails: samdarbin@ipia.sci.am, isko@ipia.sci.am\\
\end{center}

\textbf{Abstract}

Let $D$ be a strongly connected balanced bipartite directed graph of order $2a\geq 10$ other than a directed cycle. Let $x,y$ be distinct vertices in $D$. $\{x,y\}$ dominates a vertex $z$ if $x\rightarrow z$ and $y\rightarrow z$; in this case, we call the pair $\{x,y\}$ dominating.  In this paper we prove: 

  {\it If  $ max\{d(x), d(y)\}\geq 2a-2$ for every dominating pair of vertices $\{x,y\}$,  then $D$ contains cycles of all lengths $2,4, \ldots , 2a-2$ or  $D$ is isomorphic to a certain digraph of order ten which we specify.}\\


\section {Introduction} 

It is sequel to the paper \cite{[13]} by the first author. We consider digraphs (directed graphs)  in the sense of \cite{[4]}, and use standard graph theoretical terminology and notation (see Section 2 for details). A cycle passing  through all the vertices of  a digraph  is called Hamiltonian. A digraph containing a Hamiltonian cycle is called a Hamiltonian digraph.  A digraph $D$ of order $n$ is called    pancyclic if it contains cycles of every lengths $3, 4,\ldots ,  n$. Various sufficient conditions for a digraph to be Hamiltonian have been given in terms of the vertex degree of the digraph. Here we recall some of them which are  due to  Ghouila-Houri \cite{[16]}, Nash-Williams \cite{[23]}, Woodall \cite{[27]}, Meyniel \cite{[22]}, Thomassen \cite{[25]} and Darbinyan \cite{[11]}. The Meyniel theorem is a  generalization  Nash-Williams', Ghouila-Houri's and Woodall's theorems.\\
 
Bondy suggested (see \cite{[9]} by Chv\'{a}tal) the following   metaconjecture:\\

 \textbf{Metaconjecture.} {\it Almost any non-trivial condition  of a graph (digraph) which implies that the graph (digraph) is Hamiltonian  also implies that the graph (digraph) is pancyclic. (There may be  a "simple" family of exceptional graphs (digraphs)).}\\

In fact various sufficient conditions for a digraph to be Hamiltonian are also sufficient for the digraph to be pancyclic. Namely,
 in \cite{[20], [24], [10], [12]}, it was shown that if a digraph $D$ satisfies one of the above mentioned conditions for hamiltonicity digraphs, then the digraph $D$ also is pancyclic (unless some extremal cases which are characterized). For additional information on Hamiltonian and pancyclic digraphs, see, e.g, the book by Bang-Jensen and Gutin \cite{[4]} and the surveys 
\cite{[7]} by Bermond and Thomassen, \cite{[21]} by K\"{u}hn and Ostus and \cite{[17]} by Gutin.
 
 Each of aforementioned theorems  imposes a degree condition on all vertices (or, on all pairs of nonadjacent vertices). In \cite{[5]} and \cite{[3]}, it was described a type of sufficient conditions for a digraph to be Hamiltonian, in which a  degree condition requires only for some pairs of nonadjacent vertices. Let us recall  only the following theorem of them.\\

 \textbf{Theorem 1.1} (Bang-Jensen, Gutin, H.Li \cite{[5]}). {\it Let $D$ be a strongly connected digraph of order $n\geq 2$. Suppose that $min\{d(x),d(y)\}\geq n-1$ and  $d(x)+d(y)\geq 2n-1$ for any pair of non-adjacent vertices $x,y$ with a common in-neighbour. Then $D$ is Hamiltonian.}\\

A digraph $D$ is called a bipartite digraph if there exists a partition $X$, $Y$ of its vertex set into two partite sets such that every arc of $D$ has its end-vertices in different partite sets. 
It is called balanced if $|X|=|Y|$. 

There are analogies results  to the Nash-Williams, Ghouila-Houri, Woodall, Meyniel and Thomassen theorems
for balanced bipartite digraphs (see e.g., \cite{[2]} and the papers cited there).

An analogue of Theorem 1.1  for bipartite digraphs was given by R. Wang \cite{[26]} and recently a different result was given by Adamus \cite{[1]}.

\textbf{Theorem 1.2} (R. Wang \cite{[26]}). {\it Let $D$ be a strongly connected balanced bipartite digraph of order $2a$, where $a\geq 1$. Suppose that, for every dominating pair of vertices $\{x,y\}$, either $d(x)\geq 2a-1$ and $d(y)\geq a+1$ or $d(y)\geq 2a-1$ and $d(x)\geq a+1$. Then $D$ is Hamiltonian}.\\

\textbf{Theorem 1.3} (Adamus \cite{[1]}). {\it Let $D$ be a strongly connected balanced bipartite digraph of order $2a$, where $a\geq 3$. If $d(x)+d(y)\geq 3a$ for every  pair of vertices $\{x,y\}$ with a common in-neighbour or  a common out-neighbour, then $D$ is Hamiltonian}.\\

Let $D$ be a balanced bipartite digraph of order $2a$, where $a\geq 2$. For integer $k\geq 0$, we say that
 $D$ satisfies condition $B_k$ when $max\{d(x),d(y)\}\geq 2a-2+k$ for every pair of dominating vertices  $x$ and $y$.\\

Before stating the next theorems we need to define three digraphs. \\

\textbf{Example 1.} Let $D(10)$ be a bipartite digraph with partite sets $X=\{x_0,x_1,x_2,x_3,x_4\}$ and 
$Y=\{y_0,y_1,y_2,y_3,y_4\}$ satisfying the following conditions: The induced subdigraph $\langle\{x_1,x_2,x_3,y_0,y_1\}\rangle$ is a complete bipartite digraph with partite sets  $\{x_1,x_2,x_3\}$ and $\{y_0,y_1\}$;  $\{x_1,x_2,x_3\}\rightarrow\{y_2,y_3,y_4\}$; $x_4\leftrightarrow y_1$; $x_0\leftrightarrow y_0$ and $x_i\leftrightarrow y_{i+1}$ for all $i\in [1,3]$. $D(10)$ contains no other arcs.

It is easy to check that the digraph $D(10)$ is strongly connected  and satisfies condition $B_0$, but the underlying undirected graph of $D(10)$ is not 2-connected and $D(10)$ has no cycle of length 8. (It follows from the facts that $d(x_0)=d(x_4)=2$ and $x_0$ ($x_4$) is on 2-cycle). It is not difficult to check that any digraph obtained from $D(10)$ by adding a new arc whose one end-vertex is $x_0$ or $x_4$ contains no cycle of length eight. Moreover, if to $A(D)$ we add some new arcs of the type $y_ix_j$, where $i\in [2,4]$ and $j\in [1,3]$, then always we obtain a digraph which does not satisfy condition $B_0$. \\ 

\textbf{Example 2.} Let  $K_{2,3}^*$ be a complete bipartite digraph with partite sets $\{x_1,x_2\}$ and $\{y_1,y_2,y_3\}$. Let $H(8)$ be the bipartite digraph obtained from the digraph $K_{2,3}^*$ by adding  three new vertices $x_0,y_0, x_3$ and the following new arcs $x_0y_0$, $y_0x_0$, $x_0y_1$, $y_1x_0$, $x_3y_3$ and $y_3x_3$.

It is not difficult to check that the digraph $H(8)$ is strongly connected  and satisfies condition $B_0$, but the underlying undirected graph of $H(8)$ is not 2-connected  and $H(8)$ has no cycle of length 6. Hence, the bound on order of $D$ in Theorem 3.4 is sharp.\\ 

\textbf{Example 3.} Let $D(8)$ be a  bipartite digraph   with partite sets $X=\{x_0,x_1,x_2,x_3\}$ and 
$Y=\{y_0,y_1,y_2,y_3\}$, and the arc set  $A(D(8))$ contains exactly the following  arcs $y_0x_1$, $y_1x_0$, $x_2y_3$, $x_3y_2$ and  all the arcs of the following 2-cycles: 
$x_i\leftrightarrow y_i$, $i\in [0,3]$, $y_0\leftrightarrow x_2$, $y_0\leftrightarrow x_3$, $y_1\leftrightarrow x_2$ and  $y_1\leftrightarrow x_3$. 

It is easy to see that 
$$
d(x_2)=d(x_3)=d(y_0)=d(y_1)=7 \quad \hbox{and} \quad d(x_0)=d(x_1)=d(y_2)=d(y_3)=3,
$$
and  the dominating pairs in $D(8)$ are: $\{y_0,y_1\}$, $\{y_0,y_2\}$,$\{y_0,y_3\}$,$\{y_1,y_2\}$, $\{y_1,y_3\}$, $\{x_0,x_2\}$,
$\{x_0,x_3\}$,
$\{x_1,x_2\}$, $\{x_1,x_3\}$ and $\{x_2,x_3\}$. Note that every dominating pair satisfies condition $B_1$. Since $x_0y_0x_3y_2x_2$ $y_1x_0$ is a cycle of length 6 in $D(8)$, 
it is not difficult to check that $D(8)$ is strong. 

Observe that $D(8)$ is not Hamiltonian. Indeed, if $D(8)$ contains  a Hamiltonian cycle, say  $C$, then $C$ would  contain the arcs $x_1y_1$ and $x_0y_0$. Therefore, $C$ must contain the path $x_1y_1x_0y_0$ or the path $x_0y_0x_1y_1$, which is impossible since $N^-(x_0)=N^-(x_1)=\{y_0,y_1\}$.\\

For $a\geq 5$ Theorem 1.2 is an immediate consequence of the following theorem by the first author \cite{[14]}.\\

\textbf{Theorem 1.4} (Darbinyan \cite{[14]}). {\it Let $D$ be a strongly connected balanced bipartite digraph of order $2a$, where $a\geq 4$. Suppose that, for every dominating pair of vertices $\{x,y\}$, either $d(x)\geq 2a-1$  or $d(y)\geq 2a-1$. Then either $D$ is Hamiltonian or isomorphic to the digraph $D(8)$}.\\

 A balanced bipartite digraph of order $2m$ is even pancyclic  if it contains a cycle of length $2k$ 
for any $2\leq k\leq m$. 
A cycle of a balanced bipartite digraph $D$ is called pre-Hamiltonian if it contains all the vertices of $D$ except two.

Characterizations of even pancyclic bipartite tournaments was given in \cite{[6]} and \cite{[28]}. A characterization of pancyclic ordinary $k$-partite ($k\geq 3$) tournaments (respectively, pancyclic ordinary complete multipartite digraphs) was established in \cite{[18]} (respectively, in \cite{[19]}).\\

Motivated by the Bondy's metaconjecture, it is natural to set the following problem: 

\textbf{Problem.} {\it Characterize those digraphs which satisfy the conditions of Theorem 1.2 (or, 1.3 or 1.4) but are not even pancyclic.}\\
 
 In \cite{[15]}, the first author have proved the following Theorems 1.5 and 1.6.

 \textbf{Theorem 1.5} (\cite{[15]}. {\it Let $D$ be a strongly connected balanced bipartite digraph of order $2a\geq 8$ other than a directed cycle. If  $ max\{d(x), d(y)\}\geq 2a-1$ for every dominating pair of vertices $\{x,y\}$,  then  $D$ contains a cycles of all even lengths less than equal $2a$ or $D$ is isomorphic to the digraph $D(8)$.}\\

\textbf{Theorem 1.6.} (\cite{[15]}. {\it Let $D$ be a strongly connected balanced bipartite digraph of order $2a\geq 8$ which contains a cycle of length $2a-2$. If $max \{d(x), d(y)\}\geq 2a-2$ for every dominating pair of vertices $\{x,y\}$, then for any $k$, $1\leq k\leq a-1$, $D$ contains a cycle of length $2k$.}\\

In view of Theorem 1.6 it seems quite natural to ask whether a balanced bipartite digraph of order $2a$,  which satisfies condition $B_0$
  contains a pre-Hamiltonian 
cycle (i.e., a cycle of length $2a-2$). 
In this paper we prove the following theorems.\\ 

\textbf{Theorem 1.7.} {\it Let $D$ be a strongly connected balanced bipartite digraph of order $2a\geq 10$ other than the directed cycle of length $2a$. Suppose that $D$  satisfies condition $B_0$, i.e., $ max\{d(x), d(y)\}\geq 2a-2$ for every dominating pair of vertices $\{x,y\}$.  Then  $D$ contains  a cycle of lengths $2a-2$ unless $D$  is isomorphic to the digraph $D(10)$.}\\

From Theorem 1.6 and 1.7 it follows the following theorem. 

\textbf{Theorem 1.8.} {\it Let $D$ be a  balanced bipartite  digraph of order $2a\geq 10$ other than the directed cycle of length $2a$. Suppose that  $D$  satisfies condition $B_0$, i.e., $ max\{d(x), d(y)\}\geq 2a-2$ for every dominating pair of vertices $\{x,y\}$.  Then  $D$ contains  cycles of all lengths $2, 4, \ldots , 2a-2$ unless $D$ is  isomorphic to the digraph $D(10)$.}\\

\section {Terminology and Notation}

  In this paper we consider finite digraphs without loops and multiple arcs. The vertex set and the arc set of a digraph $D$ are denoted  
  by $V(D)$  and by  $A(D)$, respectively.  The order of $D$ is the number
  of its vertices. For any $x,y\in V(D)$, we also write $x\rightarrow y$, if $xy\in A(D)$. If $xy\in A(D)$, then we say that $x$ dominates $y$ or $y$ is an out-neighbour of $x$, and $x$ is an in-neighbour of $y$. 
The notation
 $x\leftrightarrow y$ denotes that $x\rightarrow y$ and 
$y\rightarrow x$ ($x\leftrightarrow y$ is called a 2-cycle). We denote by $a(x,y)$ the number of arcs with end-vertices $x$ and $y$.
For disjoint subsets $A$ and  $B$ of $V(D)$  we define $A(A\rightarrow B)$ \,
   as the set $\{xy\in A(D) | x\in A, y\in B\}$ and $A(A,B)=A(A\rightarrow B)\cup A(B\rightarrow A)$. If $x\in V(D)$
   and $A=\{x\}$ we write $x$ instead of $\{x\}$. If $A$ and $B$ are two disjoint subsets of $V(D)$ such that every
   vertex of $A$ dominates every vertex of $B$, then we say that $A$ dominates $B$, denoted by $A\rightarrow B$. 
$A\mapsto B$ means that $A\rightarrow B$ and there is no arc from $B$ to $A$.
 
We let $N^+(x)$, $N^-(x)$ denote the set of  out-neighbours, respectively the set  of in-neighbours of a vertex $x$ in a digraph $D$.  If $A\subseteq V(D)$, then $N^+(x,A)= A \cap N^+(x)$ and $N^-(x,A)=A\cap N^-(x)$. 
The out-degree of $x$ is $d^+(x)=|N^+(x)|$ and $d^-(x)=|N^-(x)|$ is the in-degree of $x$. Similarly, $d^+(x,A)=|N^+(x,A)|$ and $d^-(x,A)=|N^-(x,A)|$. The degree of a vertex $x$ in $D$ is defined as $d(x)=d^+(x)+d^-(x)$ (similarly, $d(x,A)=d^+(x,A)+d^-(x,A)$).
 
For integers $a$ and $b$, $a\leq b$, let $[a,b]$  denote the set of
all integers which are not less than $a$ and are not greater than
$b$.

The subdigraph of $D$ induced by a subset $A$ of $V(D)$ is denoted by $D\langle A\rangle$ or $\langle A\rangle$ for brevity. 
The path (respectively, the cycle) consisting of the distinct vertices $x_1,x_2,\ldots ,x_m$ ( $m\geq 2 $) and the arcs $x_ix_{i+1}$, $i\in [1,m-1]$  (respectively, $x_ix_{i+1}$, $i\in [1,m-1]$, and $x_mx_1$), is denoted by  $x_1x_2\cdots x_m$ (respectively, $x_1x_2\cdots x_mx_1$). 
We say that $x_1x_2\cdots x_m$ is a path from $x_1$ to $x_m$ or is an $(x_1,x_m)$-path.
 If $P$ is a path containing a subpath from $x$ to $y$ we let $P[x,y]$ denote that subpath.
 Similarly, if $C$ is a cycle containing vertices $x$ and $y$, $C[x,y]$ denotes the subpath of $C$ from $x$ to $y$.
 A digraph $D$ is strongly connected (or, just, strong) if there exists an $(x,y)$-path in $D$  for every ordered pair of distinct vertices $x,y$ of $D$.
  Given a vertex $x$ of a path $P$ or a cycle $C$, we denote by $x^+$ (respectively, by $x^-$) the successor (respectively, the predecessor) of $x$ (on $P$ or $C$), and in case of ambiguity, we precise $P$ or $C$ as a subscript (that is $x^+_P$ ...). 
Two distinct vertices $x$ and $y$ are adjacent if $xy\in A(D)$ or $yx\in A(D) $ (or both). 

 Let $C$ be a non-Hamiltonian cycle in a digraph $D$. An $(x,y)$-path $P$ is a $C$-bypass if $|V(P)|\geq 3$, $x\not=y$ and $V(P)\cap V(C)=\{x,y\}$. The length of the path $C[x,y]$ is the gap of $P$ with respect to $C$.

The underlying undirected graph of a digraph $D$ is the unique graph that contains an edge $xy$ if $x\to y$ or $y\to x$ (or both).

\section { Preliminaries }

 \textbf{Lemma 3.1} (Bypass Lemma 3.17, Bondy \cite{[8]}). {\it Let $D$ be a strongly connected   digraph, and let $H$ be a non-trivial proper subdigraph of $D$. If the underlying undirected graph of $D$ is  2-connected, then $D$ contains a $H$-bypass.}\\

{\it Remark.} One can prove Bypass Lemma using the proof of Theorem 5.4.2 \cite{[4]}.\\
 
\textbf{Lemma 3.2} (\cite{[14]}). {\it Let $D$ be a strongly connected balanced bipartite digraph of order $2a\geq 8$ other than the directed cycle of length $2a$. If $D$ satisfies condition $B_0$, then $D$ contains a non-Hamiltonian cycle of length at least 4. }\\

\textbf{Lemma 3.3} (\cite{[15]}). {\it Let $D$ be a strongly connected balanced bipartite digraph of order $2a\geq 8$  with partite sets $X$ and $Y$.
Assume that $D$ satisfies condition $B_0$.  Let $C=x_1y_1x_2y_2\ldots x_ky_kx_1$ be a longest non-Hamiltonian cycle in $D$, where $k\geq 2$, $x_i\in X$ and $y_i\in Y$, and $P$ be a $C$-bypass. If the gap of $P$ with respect to $C$ is equal to one, then $k=a-1$, i.e., the longest non-Hamiltonian cycle in $D$ has length  $2a-2$.}\\

\textbf{Theorem 3.4} (\cite{[13]}). {\it Let $D$ be a strongly connected balanced bipartite digraph of order $2a\geq 10$.   Assume that $D$ satisfies condition $B_0$. Then either the underlying undirected graph of $D$ is 2-connected or $D$ contains a cycle of length $2a-2$ unless $D$ is isomorphic to the digraph $D(10)$.}\\

\section { Proof of the main result}

\textbf{Proof of Theorem 1.7. } Suppose, on the contrary, that a digraph $D$ is not a directed cycle and
  satisfies the conditions of the  theorem, but contains no cycle of length $2a-2$. Let $C=x_1y_1x_2y_2\ldots x_my_mx_1$ be a longest non-Hamiltonian cycle in $D$, where $x_i\in X$ and $y_i\in Y$ (all subscripts of the vertices $x_i$ and $y_i$ are taken modulo $m$, i.e., $x_{m+i}=x_i$ and $y_{m+i}=y_i$).

 By Lemma 3.2, $D$ contains a non-Hamiltonian cycle of length at least 4, i.e., $2\leq m\leq a-2$. By Theorem 3.4, the underlying undirected graph of $D$
 is 2-connected. Therefore,
 by Bypass Lemma, $D$ contains a $C$-bypass. We choose a cycle $C$ and a $C$-bypass $P:=xu_1\ldots u_sy$ such that 

(i) $C$ is a longest non-Hamiltonian cycle in $D$;

(ii) the gap of $C$-bypass $P$ with respect to $C$ is minimum subject to (i);

(iii) the length of $P$ is minimum subject to (i) and (ii).\\
 
Without loss of generality, we assume that $x:=x_1$. Let $R:=V(D)\setminus V(C)$, $P_1:=P[u_1,u_s]$ and $C':=V(C[y_1,y^-_C])$. Note that $|R|\geq 4$. From Lemma 3.3 it follows that $|V(C[x,y])|\geq 3$.  Then $|C'|\geq s$ since $C$ is a longest non-Hamiltonian cycle in $D$. \\

We first prove that $|V(P_1)|=1$, i.e., $s=1$.

\textbf{Proof of $s=1$.}  Suppose, on the contrary, that $s\geq 2$. By the minimality of the gap $|C'|+1$, we have
$$
d(y^-_C, P_1)=d(u_s, C')=0. \eqno (1)
$$
It is not difficult to see that $s\leq 3$. Indeed, if $s\geq 4$, then from (1) it follows  that $d(y^-_C)\leq 2a-4$ and $d(u_s)\leq 2a-4$, since each of $P_1$ and $C'$ contains at least two vertices from each partite set. So, we have a contradiction, since $\{u_s,y^-_C\}\rightarrow y$.  Therefore, $s=3$ or $s=2$.

First consider the case $s=3$, i.e., $P_1=u_1u_2u_3$. From our assumption that $x=x_1$ it follows   that $u_1, u_3\in Y$ and $y\in X$. Therefore, $u_2\in X$ and $y^-_C\in Y$. Note that $x_2\in C'$ since $|C'|\geq 3$. By the minimality of the gap $|C'|+1$, we have
$$
d(x_2, \{u_1,u_3\})=d(u_2, \{y_1, y^-_C\})=0. \eqno (2)
$$
Therefore, $d(x_2)\leq 2a-4$.  From the minimality of the gap $|C'|+1$ and $P_1$ it follows that $x_1u_3\notin A(D)$.   Therefore, by (2), $d(u_3)\leq 2a-3$. This together with condition $B_0$ and $\{y^-_C, u_3\}\rightarrow y$ imply that $d(y^-_C)\geq 2a-2$. Therefore, by (2), the vertex $y^-_C$ and every vertex of $X\setminus \{u_2\}$ form a 2-cycle, i.e., $y^-_C\leftrightarrow X\setminus \{u_2\}$. In particular, $y^-_C\leftrightarrow \{x_2,z\}$, where  $z$ is an arbitrary vertex of $X\cap R\setminus \{u_2\}$. By the minimality of the gap, we have $a(z,u_3)=d^-(z, \{u_1\})=0$. Therefore, $d(z)\leq 2a-3$. This and $d(x_2)\leq 2a-4$ contradict condition $B_0$ since $\{z,x_2\}\rightarrow y^-_C$.

Now consider the case $s=2$, i.e., $P_1=u_1u_2$.
Notice that  $u_1, y\in Y$ and  $y^-_C\in X$ since $x_1\in X$. It is easy to see that $|C'|$ is even.

Assume first that $|C'|\geq 4$. By the minimality of the gap $|C'|+1$, 
$$
d(u_2,\{y_1,y^{--}_C\})=a(y_C^-,u_1)=0, \eqno (3)
$$
where $y^{--}_C$ denotes the predecessor of $y^-_C$ on $C$. 
Therefore, $d(u_2)\leq 2a-4$ and $d(y^-_C)\leq 2a-2$. Since $\{u_2,y^-_C\}\rightarrow y$, i.e., $\{u_2,y^-_C\}$ is a dominating pair, by condition $B_0$ we have,  $d(y^-_C)= 2a-2$. This together with  the second equality of (3)  imply that $y^-_C$ and every vertex of $Y\setminus \{u_1\}$ form a 2-cycle.  
In particular, $y\rightarrow y^-_C\rightarrow y_1$ and 
$y^-_C\leftrightarrow z$, where $z$ is an arbitrary vertex of $Y\cap R\setminus \{u_1\}$.
 From the minimality of the  gap $|C'|+1$ it follows that $a(z,u_2)=d^-(z,\{x_2\})=0$. Hence, $d(z)\leq 2a-3$. Now we
consider the vertex $y^{--}_C$. It is easy to see that  $y^{--}_Cy^+_C\notin A(D)$ (for otherwise, the cycle $x_1u_1u_2yy^-_Cy_1\ldots y^{--}_Cy^+_C\ldots x_1$ is longer than $C$). From this and the first equality of (3) it follows that $d(y^{--}_C)\leq 2a-3$. 
Thus, we have $\{z, y^{--}_C\}\rightarrow y^-_C$ and $max \{d(z), d(y^{--}_C)\}\leq 2a-3$, which contradict condition
 $B_0$.

Assume then that $|C'|=2$, i.e., $C'=\{y_1, x_2\}$. Then $y=y_2$. By the minimality of the  gap $|C'|+1$, 
$$
a(u_1,x_2)=a(u_2,y_1)=0,  \eqno (4)
$$
i.e., the vertices $u_1$ and $x_2$ (respectively, $u_2$ and $y_1$) are not adjacent.
Therefore,
$$
max\{d(u_1), d(x_2), d(u_2), d(y_1)\}\leq 2a-2.
$$
Without loss of generality, we may assume that $d(x_2)=2a-2$, since $\{u_2,x_2\}\rightarrow y_2$ (for otherwise,  $d(u_2)=2a-2$ and we will consider the cycle $x_1u_1u_2y_2\ldots y_mx_1$). 
Since $u_1$ and $x_2$ are non-adjacent and  $d(x_2)=2a-2$,  it follows that $x_2$ and every vertex of $Y\setminus \{u_1\}$ form a 2-cycle, i.e., $x_2\leftrightarrow Y\setminus \{u_1\}$. In particular,  $y_2\rightarrow x_2\rightarrow y_1$.
 Let $z$ be an arbitrary vertex in $ Y\cap R\setminus \{u_1\}$. By the minimality of the  gap $|C'|+1$,  $a(z, u_2)=0$ and $x_1z\notin A(D)$. Hence, $d(z)\leq 2a-3$. 
 If $y_1\rightarrow x_3$ (possibly, $x_3=x_1$), then, since $y_2\rightarrow x_2 \rightarrow y_1$, we see that  $x_1u_1u_2y_2x_2y_1x_3 \ldots x_1$ is a cycle of length $|C|+2$, a contradiction. We may therefore assume  that $y_1x_3\notin A(D)$. 
 This together with $a(y_1,u_2)=0$ (by (4)) gives $d(y_1)\leq 2a-3$. Thus we have that $\{y_1, z\}\rightarrow x_2$ and 
$max\{d(y_1), d(z)\}\leq 2a-3$, which contradict condition $B_0$. This contradiction completes the proof of $s=1$. \fbox \\\\

From $s=1$ it follows that $u_1\in Y$ and $y\in X$ since $x_1\in X$. Without loss of generality, we may assume that $y:=x_r$. From now on, let $y:=u_1$.\\

 Now we divide the proof of the theorem into two parts: $|C'|=1$ and $|C'|\geq2$.\\

\textbf{Part I}. $|C'|=1$, i.e., $r=2$ and $x_1\rightarrow y\rightarrow x_2$.

By condition $B_0$, 
$max\{d(y), d(y_1)\}\geq 2a-2$ since $\{y,y_1\}\rightarrow x_2$. Without loss of generality, assume that $d(y)\geq 2a-2$. For this part we first  prove Claims 1-5 below.\\

\textbf{Claim 1.} If $x\in R\cap X$ and $x\leftrightarrow y$, then $d(x)\leq 2a-3$ and $d(x_1)\geq 2a-2$.

\textbf{Proof of Claim 1.} Assume that $x\in R\cap X$ and $x\leftrightarrow y$, but $d(x)\geq 2a-2$. 
It is easy to see that $a(x,y_1)=0$ since $C$ is a longest non-Hamiltonian cycle in $D$. This and  $d(x)\geq 2a-2$
imply that the vertex $x$ and every vertex of $Y\setminus \{y_1\}$ form a 2-cycle. 
In particular, $x\leftrightarrow \{y_0,y_2,y_m\}$ (possibly, $y_2=y_m$), where $y_0$ is an arbitrary vertex of $Y\cap R\setminus \{y\}$.
Using this, it is easy to check that 
$$
a(x_1,y_0)= d^-(x_1,\{y\})=d^-(y,\{x_2\})=0
$$
since $C$ is a longest non-Hamiltonian cycle in $D$. From the last equalities we have $d(x_1)\leq 2a-3$, and $y\leftrightarrow x_0$, where $x_0$ is an arbitrary vertex of $X\cap R\setminus \{x\}$, since $d(y)\geq 2a-2$. Since $C$ is a longest non-Hamiltonian cycle in $D$ and since $x\leftrightarrow \{y_0,y_2,y_m\}$, it follows 
that $d(x_0,\{y_1,y_0\})=0$. Therefore, $d(x_0)\leq 2a-4$, which contradicts condition $B_0$ since 
$max\{d(x_1), d(x_0)\}\leq 2a-3$ and $\{x_0,x_1\}\rightarrow y$. This contradiction proves that $d(x)\leq 2a-3$. From this and condition 
$B_0$ it follows that $d(x_1)\geq 2a-2$ since $\{x,x_1\}\rightarrow y$. Claim 1 is proved. \fbox \\\\

 Claim 2 follows immediately from Claim 1 and condition $B_0$.\\

\textbf{Claim 2.} There are no two distinct vertices  $x, x_0\in R\cap X$ such that  $x\leftrightarrow y$ and 
$x_0\leftrightarrow y$. \fbox \\\\

\textbf{Claim 3.} If $y\leftrightarrow x$ for some $x\in R\cap X$, then $a(y,z)=1$ 
for all
 $z\in R\cap X\setminus \{x\}$.

\textbf{Proof of Claim 3.} Let $x\leftrightarrow y$ for some  $x\in R\cap X$. Then Claim 2 implies that $a(y,z)\leq 1$ for all
 $z\in R\cap X\setminus \{x\}$. 

We want to show that $a(y,z)= 1$ for all
 $z\in R\cap X\setminus \{x\}$. Assume that this is not the case. Then $a(y,x_0)=0$ for some 
 $x_0\in R\cap X\setminus \{x\}$. This together with Claims 1 and 2 thus imply that $|R|=4$ since $d(y)\geq 2a-2$ by our assumption.  
Then $R\cap X= \{x,x_0\}$. 
Let $R\cap Y= \{y,y_0\}$.
Since $a(y,x_0)=0$ and $d(y)\geq 2a-2$, it follows that 
$$
y\leftrightarrow \{x_1,x_2,\ldots , x_m\}.  \eqno (5)
$$
 Since $C$ is a longest non-Hamiltonian cycle in $D$, (5) implies that
$d(x, \{y_1,y_2,\ldots , y_m\})=0$,    
which in turn implies that the vertices $x$ and $y_0$ are adjacent and $d(x)\leq 2a-4$ since $UG(D)$ is 2-connected
and $m\geq 2$.  By the assumption of Claim 3 and (5), $\{x,x_i\}\rightarrow y$ for all $i\in [1,m]$. Now using 
condition $B_0$,  we obtain
$$
d(x_i)\geq 2a-2,  \quad \hbox{for all} \quad i\in [1,m]. \eqno (6)
$$
The remainder of the proof of Claim 3 is divided into two subcases depending on the value  of $a(x,y_0)$.

\textbf{Case 1.} $a(x,y_0)=2$, i.e., $x\leftrightarrow y_0$.

Since $C$ is a longest non-Hamiltonian cycle in $D$, from  (5) and $y\leftrightarrow x$, $x\leftrightarrow y_0$ it follows that 
$$
d(y_0, \{x_1,x_2,\ldots , x_m\})=0.     \eqno (7)
$$
Therefore, since $UG(D)$ is 2-connected, the vertices $y_0$ and $x_0$ are adjacent. Using (5), it is not difficult to see that:

If $x_0\rightarrow y_0$, then $d^-(x_0, \{y_1,y_2,\ldots , y_m\})=0$, and
if $y_0\rightarrow x_0$, then $d^+(x_0, \{y_1,y_2,\ldots , y_m\})=0$.

In both cases we have $a(y_i,x_0)\leq 1$ for all $i\in [1,m]$.  Together with $a(y_i,x)=0$ this implies that $d(y_i)\leq 2a-3$. On the other hand,
from (6) and (7) it follows that every vertex $x_i$, $i\in [1,m]$, and every vertex of $Y\setminus \{y_0\}$ form a 2-cycle. 
In particular, $\{y_1,y_2\}\rightarrow x_2$, which contradicts condition $B_0$, since $max\{d(y_1),d(y_2)\}\leq 2a-3$.

\textbf{Case 2.} $a(x,y_0)=1$, i.e., $y_0\mapsto x$ or $x\mapsto y_0$.

Let  $y_0\mapsto x$.  Since $C$ is a longest non-Hamiltonian cycle in $D$ and  (5), we have 
$$
d^-(y_0, \{x_1,x_2,\ldots , x_m\})=0.
$$
Hence $x_0\rightarrow y_0$ since $D$ is strong and  $y_0\mapsto x$. Now using (5), we obtain  
$$
d^-(x_0, \{y_1,y_2,\ldots , y_m\})=0.
$$
The last two equalities together with $a(y,x_0)=0$ and $y_0\mapsto x$ imply that 
$$A(V(C)\cup \{x,y\}\rightarrow \{x_0, y_0\})=\emptyset ,$$ 
which contradicts that $D$ is strong.

Let now $x\mapsto y_0$. Again using (5), it is easy to see that 
$d^+(y_0, \{x,x_1,x_2,\ldots , x_m\})=0$. Therefore, $y_0\rightarrow x_0$ since $D$ is strong.  Together with (5) this implies that  
$d^+(x_0, \{y_1,y_2,\ldots , y_m,y\})=0$.
Therefore $x_0\rightarrow y_0$ since $D$ is strong, and $d(x_0)\leq 2a-4$ since $m\geq 2$. Thus,  $d(x)\leq 2a-3$ (Claim 1) and $d(x_0)\leq 2a-4$, which contradicts condition $B_0$ since $\{x,x_0\}\rightarrow y_0$. This contradiction completes the proof of Claim 3. \fbox \\\\

\textbf{Claim 4.} If $y\leftrightarrow x$ for some $x\in R\cap X$, then $d^-(y,R\cap X\setminus \{x\})=0.$ 

\textbf{Proof of Claim 4.} Assume that  the claim is not true, i.e., there exist vertices  $x\in R\cap X$ and  $u\in R\cap X\setminus \{x\}$ such that $y\leftrightarrow x$ and
$u\rightarrow y$. Then, by Claim 2,  $yu\notin A(D)$. Notice that $\{x,u\}\rightarrow y$. Since $d(x)\leq 2a-3$  (Claim 1), condition $B_0$ implies that  $d(u)\geq 2a-2$. It is clear that $y_1u\notin A(D)$ (if $y_1\rightarrow u$, then $x_1y_1uyx_2y_3\ldots y_mx_1$ is a cycle of length $2m+2$, a contradiction). Using this, $yu\notin A(D)$ and $d(u)\geq 2a-2$ we conclude that $uy_1\in A(D)$, $d(u)=2a-2$ and the vertex $u$ together with every vertex of $Y\setminus \{y,y_1\}$ forms a 2-cycle.
In particular, $u\leftrightarrow \{y_2,y_m, v\}$, where $v\in Y\cap R\setminus \{y\}$ (possibly, $y_2=y_m$). 
Using this, 
it is not difficult to show that $a(x_1,v)=0$ and $y x_1\notin A(D)$. 
Indeed, if $x_1\rightarrow v$, then $x_1vuy_1x_2\ldots y_mx_1$ is a cycle of length $2m+2$; 
if $v\rightarrow x_1$, then
 $y_muvx_1y_1\ldots $ $x_my_m$ is a cycle of length $2m+2$; 
if $y\rightarrow x_1$, then $y_muyx_1y_1\ldots x_my_m$ is a cycle of length $2m+2$. In each case we obtain a contradiction. Hence, $a(x_1,v)=0$ and $yx_1\notin A(D)$. From this it follows that $d(x_1)\leq 2a-3$, which contradicts that $d(x_1)\geq 2a-2$ (Claim 1). Claim 4 is proved.   \fbox \\\\

\textbf{Claim 5.} There is no $x\in X\cap R$ such that $y\leftrightarrow x$, i.e., in subdigraph $D\langle R\rangle $ through  the vertex $y$ there is no cycle of length two. 

\textbf{Proof of Claim 5.} Assume that  the claim is not true, i.e., there exists a vertex $x\in X\cap R$ such that 
$y\leftrightarrow x$.
By Claims 1, 3 and 4 we have
$$
d(x)\leq 2a-3, \quad d(x_1)\geq 2a-2 \quad \hbox{and} \quad y\mapsto X\cap R \setminus \{x\}.    \eqno (8)
$$

\textbf{Case 1}. $x_2y\notin D$.

Then from the last expression of (8) and $d(y)\geq 2a-2$ we conclude that $|R|=4$, i.e., $m=a-2\geq 3$, and the vertex $y$ and every vertex of $\{x_1,x_2,x_3,x_4,\ldots , x_m\}\setminus \{x_2\}$ form a 2-cycle, i.e., 
$$
 y\leftrightarrow \{x_1,x_3,x_4,\ldots , x_m\} \quad \hbox{and} \quad y\mapsto x_2. \eqno (9)
$$

Put $R=\{y,x,y_0, x_0\}$, where $x_0\in X$ and $y_0\in Y$. From (8) it follows that 
$y\mapsto x_0$. Now using (9), it is not difficult to see that 
$$
d^-(x,\{y_1,y_2,\ldots ,y_m\})=d^+(x,\{y_1,y_2,\ldots ,y_m\}-\{y_2\})=0. \eqno (10)
$$
Indeed, if $y_i\rightarrow x$ and $i\in [1,m]$ (respectively, $x\rightarrow y_j$ and $j\in [1,m]\setminus \{2\}$), then
$x_iy_ixyx_{i+1}\ldots x_i$ (respectively, $x_jyxy_j\ldots y_{j-1}x_j$)
is a cycle of length $2a-2$, a contradiction.

Similarly,
$$
d^+(x_0,\{y_1,y_2,\ldots ,y_m\}-\{y_2\})=0. \eqno (11)
 $$ 
In particular, (10) implies that
 $$
d(x,\{y_1,y_2,\ldots ,y_m\}-\{y_2\})=0, 
 $$
i.e., the vertices $x$ and $y_i$, $i\notin \{0,2\}$, are not adjacent. From this, (11) and $y\mapsto x_0$ it follows that 
$$
max\{d(x_0),d(x),d(y_j)\}\leq 2a-3, \quad \hbox{for all}\quad j\notin \{0,2\}. \eqno (12)
$$
Using (12) and condition $B_0$, we obtain that for all $u\in X$ and $v\in Y$ the following holds
$$
d^-(u,\{y_1,y_m\})\leq 1  \quad \hbox{and}\quad d^-(v,\{x_0,x\})\leq 1. \eqno (13)
$$
Now we divide this case into four subcases.

\textbf{Subcase 1.1.} $x\leftrightarrow y_0$.

Since $x\leftrightarrow y$ and $x\leftrightarrow y_0$, using (9) and the fact that $m\geq 3$, it is not difficult to check that the vertices $x_1$ and $y_0$ are not adjacent (for otherwise, $D$ would contain a cycle of length $2a-2$, a contradiction).  Together with the first inequality of (13) (when $u=x_1$) this implies that $d(x_1)\leq 2a-3$, which contradicts Claim 1.

\textbf{Subcase 1.2.} $x\mapsto y_0$.

Then, by the second inequality of (13) when $v=y_0$, $x_0y_0\notin A(D)$. This together with $x_0y\notin A(D)$ and (11) imply that $x_0\rightarrow y_2$ since $D$ is strong. It is easy to see that $y_0x_2\notin A(D)$ (for otherwise, $x_1yxy_0x_2y_2\ldots y_mx_1$ is a cycle of length
$2a-2$, a contradiction). Combining this with $x_2y\notin A(D)$ and 
$d^-(x_2,\{y_1,y_m\})\leq 1$ (by (13)) we obtain $d(x_2)\leq 2a-3$. Thus, $d(x_2)\leq 2a-3$ and $d(x_0)\leq 2a-3$ (by (12))
and $\{x_0,x_2\}\rightarrow y_2$, which contradict condition $B_0$.

\textbf{Subcase 1.3.} $y_0\mapsto x$.

If $x_i\rightarrow y_0$ and $i\in [1,m]$, then using (9) we obtain that the cycle $x_iy_0xyx_{i+1}\ldots x_i$ has length $2a-2$, which is a contradiction. We may therefore assume  that $d^-(y_0,\{x,x_1,x_2,\ldots ,x_m\})=0$. Hence, $x_0\rightarrow y_0$
since $D$ is strong and  $y_0\mapsto x$. 
Since $x_m\rightarrow y\rightarrow x_2$ and $C$ is a longest non-Hamiltonian cycle in $D$, it follows that the vertices $x_1$ and $y_0$ are not adjacent. Combing this with $d^-(x_1,\{y_1,y_m\})\leq 1$ (by (13)) we obtain that $d(x_1)\leq 2a-3$, a contradiction to Claim 1.

\textbf{Subcase 1.4.} The vertices $x$ and $y_0$ are not adjacent.

Since the underlying undirected graph of $D$ is 2-connected, from $a(x,y_0)=0$ and (10) it follows that $x\rightarrow y_2$. Together with (13) this imply that $x_0y_2\notin A(D)$. Therefore, by (11), we have that
$d^+(x_0,\{y,y_1,y_2,\ldots , y_m\})=0$. Hence, $x_0\rightarrow y_0$ since $D$ is strong. If $y_0\rightarrow x_2$, then $x_1yx_0y_0 x_2y_2\ldots y_mx_1$ is a cycle of length $2a-2$, which is a contradiction. We may therefore assume  that $y_0x_2\notin A(D)$. Combining this with $x_2y\notin A(D)$ and (13) we obtain that $d(x_2)\leq 2a-3$, which contradicts condition $B_0$ since $\{x_2,x\}\rightarrow y_2$ and $d(x)\leq 2a-3$ (Claim 1). 
 The discussion of Case 1 is completed.

\textbf{Case 2.} $x_2\rightarrow y$.

Since $y\mapsto X\cap R\setminus \{x\}$ (by (8)) and  $d(y)\geq 2a-2$, it follows that the vertices $y$ and $x_i$, where 
$i\in [1,m]$, are adjacent and $y\rightarrow \{x_1,x_2, \ldots , x_m\}$ or $\{x_1,x_2, \ldots , x_m\}\rightarrow y$.
Therefore, without loss of generality, we may assume that 

$$x_1\rightarrow y  \quad \hbox{and} \quad x_i\leftrightarrow y \quad \hbox{for all} \quad i\in [2,m]   \eqno (14)
$$ 
(for otherwise we will have the considered Case 1).  Let $x_0$ be an arbitrary vertex in $X\cap R\setminus \{x\}$. Then by (8) we have $y\mapsto x_0$. It is clear that
$$
 d^+(x_0,\{y,y_1,y_2,\ldots ,y_m\})=d(x,\{y_1,y_2,\ldots ,y_{m-1}\})=0 \quad \hbox{and} \quad x y_m\notin A(D).  \eqno (15)
$$
 This implies that $d(x_0)\leq 2a-3$. 
Since $D$ is strong and  (15), there is a vertex $y_0\in Y\cap R\setminus \{y\}$   such that $x_0\rightarrow y_0$. 
Now, since $y\rightarrow x_0\rightarrow y_0$ and $x_i\rightarrow y$ for all $i\in [1,m]$ (by (14)), we conclude that 
$d^+(y_0,\{x_1,x_2, \ldots ,x_m\})=0$ (for otherwise, for some $i\in [1,m]$, $y_0\rightarrow x_i$ and 
$x_{i-1}yx_0y_0x_i\ldots y_{i-2}x_{i-1}$ is a cycle of length $2m+2$, a contradiction). 
The last equality together with 
$d^+(x_0,\{y,y_1,y_2,\ldots ,y_m\})=0$ (by (15))
imply that $y_0\rightarrow x$ since $D$ is strong (for otherwise,  
$A(\{x_0,y_0\}\rightarrow V(D)\setminus \{x_0,y_0\})=\emptyset$, which contradicts that $D$ is strong). 
Then using the facts that
 $y\rightarrow x_0\rightarrow y_0$ and $y_0\rightarrow x\rightarrow y$,
it is not difficult to show that $x_1$ and $y_0$ are not adjacent. 
Indeed, by (14) we have that: If $x_1\rightarrow y_0$, then 
$x_{1}y_0xyx_2y_2\ldots  y_{m}x_{1}$ is a cycle of length $2m+2$; and if $y_0\rightarrow x_1$, then $x_{m}yx_0y_0x_1y_1\ldots $ $y_{m-1}x_{m}$ is a cycle of length $2m+2$, in both cases we have a contradiction. So, $x_1$ and $y_0$ are not adjacent.
Together with $d(x_1)\geq 2a-2$ (Claim 1) this implies that $y_1\rightarrow x_1$. On the other hand, since $d(x,\{y_1,y_m\})=d^+(x_0,\{y_1,y_m\})=0$, we have that
$max\{d(y_1),d(y_m)\}\leq 2a-3$, which contradicts condition $B_0$, because of $\{y_m,y_1\}\rightarrow x_1$. Claim 5 is proved. \fbox \\\\

Now we can finish the discussion of Part I.\\

 From Claim 5 it follows that in $D\langle R\rangle$ there is no cycle of length two through the vertex $y$. Then, since $d(y)\geq 2a-2$ and $|R|\geq 4$, it follows that $|R|=4$.   
 Put $X\cap R=\{x,x_0\}$ and $Y\cap R=\{y,y_0\}$. Then $a(y,x)=a(y,x_0)=1$ and the vertex $y$ and every vertex of $X\cap C$ form a 2-cycle, i.e., 
$$
y\leftrightarrow \{x_1,x_2,\ldots , x_m\}. \eqno (16)
$$

 First consider the case $d^+(y,\{x,x_0\})\geq 1$. Assume, without loss of generality, that $y\mapsto x$. Then, by (16), 
$d^+(x,\{y_1,y_2,\ldots ,y_m\})=0$.
Together with $xy\notin A(D)$ this implies that $x\rightarrow y_0$ since $D$ is strong. 
Therefore
$d(x)\leq 2a-3$ since $|Y\cap C|\geq 2$.
 By (16), it is clear that 
$d^+(y_0,\{x_1,x_2, \ldots ,x_m\})=0$. 
If $y\rightarrow x_0$, then analogously we obtain  that 
$d^+(x_0,\{y,y_1,y_2,\ldots ,y_m\})=0$ and $x_0\rightarrow y_0$, $d(x_0)\leq 2a-3$, which contradicts condition $B_0$ since
$max\{d(x),d(x_0)\}\leq 2a-3$ and $\{x,x_0\}\rightarrow y_0$.
We may assume therefore that $yx_0\notin A(D)$. 
Then $x_0\rightarrow y$ (by $a(y,x_0)=1$), 
$d^-(x_0,\{y,y_1,y_2,\ldots ,y_m\})=0$ (by (16)) and hence, $y_0\rightarrow x_0$ since $D$ is strong. Now it  is not difficult to show that
$$
d^-(y_0,\{x_1,x_2,\ldots ,x_m\})=d^-(x_0,\{y_1,y_2,\ldots ,y_m\})=d^+(x,\{y_1,y_2,\ldots ,y_m\})=0.
$$
Therefore $d(x_0)\leq 2a-3$ and $d(y_i)\leq 2a-3$. Since for all $i\in [1,m]$, $\{x_i,x_0\}\rightarrow y$ and $d(x_0)\leq 2a-3$, 
from condition $B_0$ it follows that $d(x_i)\geq 2a-2$ for all $i\in [1,m]$. This together with $a(x_i,y_0)=0$ (by (16) and the last equalities) imply that $y_i\leftrightarrow x_i$. 
Thus, $\{y_{i-1},y_i\}\rightarrow x_i$ and 
$max\{d(y_{i-1}),d(y_i)\}\leq 2a-3$, which is a contradiction.

Now consider the case $d^+(y,\{x,x_0\})=0$. Then $\{x,x_0\}\rightarrow y$, because of $a(y,x)=a(y,x_0)=1$, i.e., $\{x,x_0\}$ is a dominating pair. It is clear that
$d^-(x,\{y_1,y_2\})=d^-(x_0,\{y_1,y_2\})=0$. This together with $d^+(y,\{x_0,x\})=0$ imply that 
 $max\{d(x),d(x_0)\}\leq 2a-3$, which is a contradiction because of $\{x,x_0\}\rightarrow y$. This completes the discussion of the part $|C'|=1$.\\

\textbf{Part 2}. $|C'|\geq 2$.

Then $|C'|\geq 3$ since $|C'|$ is odd.  
For this part we first will prove Claims 6-8.\\

\textbf{Claim 6.} If $|C'|\geq 3$, then the following holds:

(a). $d(y)\leq 2a-3$ and $d(y_{r-1})\geq 2a-2$;

(b). There is no  $x\in X\cap R$ such that $x\leftrightarrow y_{r-1}$, 
i.e., $a(y_{r-1},x)\leq 1$ for all $x\in X\cap R$;

(c). $a(y_{r-1},x)=1$ for all $x\in X\cap R$, $|R|=4$, $d(y_{r-1})=2a-2$ and the vertex $y_{r-1}$ together with every vertex of $X\cap V(C)$ forms a 2-cycle. In particular, $x_r\leftrightarrow y_{r-1}$ and $y_{r-1}\leftrightarrow x_2$;

(d). $d^-(y_{r-1},\{x,x_0\})=0$ and $y_{r-1}\mapsto \{x,x_0\}$;

(e). $max\{d(x_0), d(x)\}\leq 2a-3$ and $d^-(v,\{x_0,x\})\leq 1$ for all $v\in Y$;  

(f). $|C'|=3$, i.e., $r=3$.

\textbf{Proof of Claim 6.}

 \textbf{ (a).} Suppose on the contrary, that $d(y)\geq 2a-2$. Then, since $d(y,C')=0$, we have that $r=3$ and the vertex $y$ together with every vertex of $X\setminus \{x_2\}$ forms a 2-cycle. In particular, $y\leftrightarrow x_1$ and $y\leftrightarrow x_3$. 
It follows that for some $i\in [3,m]$, $x_i\rightarrow y\rightarrow x_{i+1}$, which contradicts that $C$-bypass $P$ has the minimum gap among the gaps of all $C$-bypasses. Therefore, $d(y)\leq 2a-3$. Now, since $\{y,y_{r-1}\}\rightarrow x_r$, from condition $B_0$ it follows that
$d(y_{r-1})\geq 2a-2$. 
   \fbox \\\\

\textbf{(b).} 
Suppose that  Claim 6(b) is falls, i.e., there is a vertex $x\in X\cap R$ such that $x\leftrightarrow y_{r-1}$. From the minimality of the gap $|C'|+1$ it follows  that
$$
a(x,y)=d^-(x,\{y_1,y_2, \ldots , y_{r-2}\})=0 \quad \hbox{and} \quad xy_r\notin A(D).    \eqno (17)
$$
From this we have  $d(x)\leq 2a-4$. This together with  condition $B_0$  and $\{x,x_{r-1}\}\rightarrow y_{r-1}$ imply that $d(x_{r-1})\geq 2a-2$.
Therefore, since $a(x_{r-1},y)=0$, it follows that  

(i) {\it the vertex $x_{r-1}$ and every vertex of $Y\setminus \{y\}$ form a 2-cycle.}\\ 

Let $y'$ be an arbitrary vertex in $Y\cap R\setminus \{y\}$. By (i), $x_{r-1}\leftrightarrow y'$.
By the minimality of the gap $|C'|+1$, $d^+(y', \{x,x_r\})= d^-(y', \{x_1\})=0$.
 Therefore, $d(y')\leq 2a-3$. Since $\{y',y_{r-2}\}\rightarrow x_{r-1}$, condition $B_0$ implies that 
$$ 
d(y_{r-2})\geq 2a-2.    \eqno (18)
$$
 First consider the case  $r\geq 4$. Then $x_{r-2}\in C'$ and it is not difficult to see that
 $x_{r-2} y'\notin A(D)$ and  $y_{r-2}x\notin A(D)$. Now, since  $a(x_{r-2}, y)=0$, we have  
 $d(x_{r-2})\leq 2a-3$. Using condition $B_0$ and the fact that $d(x)\leq 2a-4$, we conclude that $x y_{r-2}\notin A(D)$. 
Therefore, by (17), $x$ and $y_{r-2}$ are not adjacent.
This together with $d(y_{r-2})\geq 2a-2$ (by (18)) imply that $y_{r-2}\leftrightarrow x_0$, where $x_0$ is an arbitrary vertex in   $R\cap X\setminus \{x\}$. By the minimality of the gap $|C'|+1$, we have $x_0y_{r-1}\notin A(D)$, $y_{r-3}x_0\notin A(D)$ and $x_0 y'\notin A(D)$.
Therefore, $d(x_0)\leq 2a-3$, which contradicts condition $B_0$, since $\{x_0,x_{r-2}\}\rightarrow y_{r-2}$ and
$d(x_{r-2})\leq 2a-3$.

Now consider the case $r=3$. By (i), $x_2\leftrightarrow y_3$ and $x_2\leftrightarrow y_2$.  If
 $y_2\rightarrow x_4$ (possibly, $x_4=x_1$) (respectively,  $x_3\rightarrow y_2$), then 
the cycle
 $Q:=x_1yx_3y_3x_2y_2x_4\ldots y_mx_1$ (respectively,  $Q:=x_1yx_3y_2x_2y_3x_4\ldots y_mx_1 $)  has length $2m$, the vertex $y_1$ is not on this cycle and  $x_1\rightarrow y_1\rightarrow x_2$ is a $Q$-bypass whose gap with respect to $Q$ is equal to 4 and $d(y_1)\geq 2a-2$, (by (18) since $r-2=1$), which contradicts Claim 6(a). 
We may therefore assume that  $y_2 x_4\notin A(D)$ and $x_3y_2\notin A(D)$.   
Combining this with $d(y_2)\geq 2a-2$ (by Claim 6(a) and $r=3$) we obtain that $d(y_2)= 2a-2$ and $y_2$ together with every vertex of $X\setminus \{x_3,x_4\}$ forms a 2-cycle. In particular, 
 $y_2\leftrightarrow x_0$. 
On the other hand, from the minimality of the gap $|C'|+1$ it follows that $x_0, y$ are not adjacent and $y_1 x_0\notin A(D)$. Therefore,
 $d(x_0)\leq 2a-3$, which is a contradiction, since $d(x)\leq 2a-3$ and $\{x,x_0\}\rightarrow y_2$.  \fbox\\\\

\textbf{(c)}. Claim 6(c) is an immediate corollary of Claims 6(a) and 6(b).   \fbox \\\\

From now on, we assume that $X\cap R =\{x,x_0\}$ and $Y\cap R =\{y,y_0\}$.\\

\textbf{(d).} Suppose, on the contrary, that there exists a vertex in $\{x,x_0\}$, say $x$, such that $x\rightarrow y_{r-1}$.
From the minimality  of the  gap $|C'|+1$ it follows that $d^-(x,\{y,y_1\})=0$. This together with $y_{r-1} x\notin A(D)$ (Claim 6(b)) imply that $d(x)\leq 2a-3$. 
Therefore, using Claim 6(c) and  condition $B_0$, we obtain that $d(x_2)\geq 2a-2$ and $d(x_{r-1})\geq 2a-2$ (possibly, $x_2=x_{r-1}$) since $\{x,x_2,x_{r-1}\}\rightarrow y_{r-1}$. Together with $d(y,\{x_2,x_{r-1}\})=0$ this implies that $x_{r-1}\rightarrow y_r$ and 
$x_2\leftrightarrow y_0$.  
Because of gap minimality, 
 $d^-(y_0,\{x_1\})=d^+(y_0,\{x,x_r\})=0$. Therefore, $d(y_0)\leq 2a-3$, and hence, by condition $B_0$, $d(y_1)\geq 2a-2$ since $\{y_1,y_0\}\rightarrow x_2$. 
Since $x_r\rightarrow y_{r-1}$ and $y_{r-1}\rightarrow x_2$ (Claim 6(c)) and $x_{r-1}\rightarrow y_r$, the cycle $Q:=x_1yx_ry_{r-1}x_2\ldots x_{r-1}y_r\ldots y_mx_1$ has length equal to $2a-4$.
Observe that $y_1\notin V(Q)$ and $x_1\rightarrow y_1\rightarrow x_2$ is a $Q$-bypass, whose gap with respect to $Q$ is   equal to 4, but $d(y_1)\geq 2a-2$, which contradicts Claim 6(a).  \fbox \\\\

\textbf{(e).} By Claim 6(d),  $y_{r-1}\mapsto \{x,x_0\}$. Therefore, because of gap minimality, we have  
$$
A(\{x,x_0\}\rightarrow \{y,y_{r-1},y_r\})=\emptyset.
$$
 Therefore, $max\{d(x_0), d(x)\}\leq 2a-3$ and, by condition $B_0$,
the vertices $x$, $x_0$ does not form a dominating pair, i.e., $d^-(v,\{x_0,x\})\leq 1$ for all $v\in Y$.   \fbox \\\\ 
 
\textbf{(f).} Suppose, on the contrary, that is $|C'|>3$. Then $|C'|\geq 5$, i.e., $r\geq 4$, since $|C'|$ is odd.
By Claim 6(c), $x_i\leftrightarrow y_{r-1}$ for all $i\in [1,m]$. This together with condition $B_0$ imply that  
$d(x_i)\geq 2a-2$ for all $i\in [1,m]$ maybe except one. In particular,
 $d(x_{r-1})\geq 2a-2$ or $d(x_{r-2})\geq 2a-2$.

Fist consider the case  $d(x_{r-1})\geq 2a-2$. From this and $a(y,x_{r-1})=0$ it follows that $x_{r-1}\leftrightarrow y_{r}$.
Using this and the fact that
$x_r\rightarrow y_{r-1}\rightarrow x_2$ (Claim 6(c)), we see that the  cycle 
$Q:=x_1yx_ry_{r-1}x_2y_2\ldots x_{r-1}y_r$ $x_{r+1}  \ldots x_1$ has length equal to $2a-4$. Observe that $x_1\rightarrow y_{1}\rightarrow x_2$ is a $Q$-bypass whose  gap with respect to $Q$ is equal to 4, which  contradicts the choice of the cycle $C$ and $C$-bypass $P$.

Now consider the case  $d(x_{r-1})\leq 2a-3$. Then $d(x_{r-2})\geq 2a-2$. Observe that the vertex $x_{r-2}$ and every vertex of $Y$ other than $y$ form a 2-cycle (since $x_{r-2}$ and $y$ are not adjacent), in particular, $y_0\leftrightarrow x_{r-2}$. It is easy to see that $x_{r-2}=x_2$, i.e., $r=4$. Indeed, if $r-2\geq 3$, then $x_{r-3}\in C'$ and $x_{r-3}y_0\notin A(D)$ because of the minimality  of the gap $|C'|+1$ and $y_0\leftrightarrow x_{r-2}$. 
This together with $a(y,x_{r-3})=0$ gives $d(x_{r-3})\leq 2a-3$. Thus we have that the vertices $x_{r-3}$ and $x_{r-1}$ both have degree less than $2a-2$, which contradicts the fact that at most one vertex $x_i$, $i\in [1,m]$ maybe has degree less that $2a-2$.  
 Thus, $r=4$. By the above observations,  $y_4\rightarrow x_2$ and $y_3\rightarrow x_5$. 
Therefore, $Q:=x_1yx_4y_4x_2y_2x_3 y_3x_5\ldots x_1$ is a cycle of length $2a-4$. Notice that $y_1\notin V(Q)$ and the $Q$-bypass $x_1\rightarrow y_1\rightarrow x_2$ has gap with respect to $Q$  equal to 4. This  contradicts the choice of $C$ and  $C$-bypass $P$.  \fbox \\\\

From Claims  6(c), 6(d) and 6(f)1 it follows that 

(ii) {\it  If $|C'|\geq 3$, then $|C'|= 3$, $|R|=4$, $y_2\mapsto \{x,x_0\}$, $d(y)\leq 2a-3$, 
 $d(y_2)= 2a-2$ and the vertex $y_2$ together with every vertex of $X\cap V(C)$ forms a 2-cycle.}\\

\textbf{Claim 7.} If $|C'|=3$, then $d(y_0)\leq 2a-3$, (recall that $\{y_0\}=Y\cap R\setminus \{y\}$).

\textbf{Proof of Claim 7.} By Claim 6((e),  $d^-(y_0,\{x,x_0\})\leq 1$. Therefore, if $x_2$ and $y_0$ are not adjacent, then 
$d(y_0)\leq 2a-3$. We may therefore assume that $x_2$ and $y_0$ are  adjacent. Then $y_0\rightarrow x_2$ or $x_2\rightarrow y_0$.
It is easy to see that if  $y_0\rightarrow x_2$, then  $x_1y_0\notin A(D)$; and if  $x_2\rightarrow y_0$, then $y_0x_3\notin A(D)$. Therefore, $d(y_0,\{x_1,x_2,x_3\})\leq 4$. This and $d^-(y_0,\{x_0,x\})\leq 1$  imply that 
$d(y_0)\leq 2a-3$. \fbox \\\\

Combining Claims 6(a), 6(e) and 7 we obtain that if $|C'|=3$, then 
$$ 
max\{d(y),d(y_0), d(x), d(x_0)\}\leq 2a-3,    \eqno (19)
$$
 in particular, by condition $B_0$, we have
$$
max\{ d^-(u,\{y,y_0\}), d^-(v,\{x,x_0\})\}\leq 1 \quad \hbox{for all} \quad u\in X \quad \hbox{and} \quad v\in Y. \eqno (20)
$$

\textbf{Claim 8.}  If $|C'|=3$, then the following holds:

(a). $d(x_2)\leq 2a-3$ and $d(x_3)\geq 2a-2$;

(b). $d(y_1)\geq 2a-2$;

(c). If $x\rightarrow y_k$ or $x_0\rightarrow y_k$, where $k\in [3,m]$, 
then $x_k y_1\notin A(D)$.

\textbf{Proof of Claim 8.} 

\textbf{(a).}
Suppose, on the contrary, that $|C'|=3$ and $d(x_2)\geq 2a-2$. This and $a(x_2,y)=0$ imply that  $d(x_2)= 2a-2$ and $x_2$ together with  every vertex of $Y\setminus \{y\}$ form a 2-cycle. In particular, $x_2\leftrightarrow y_0$ and 
 $x_2\leftrightarrow y_3$. Since  $d(y_0)\leq 2a-3$ (by (19)) and $\{y_1,y_0\}\rightarrow x_2$, from condition $B_0$ it follows that $d(y_1)\geq 2a-2$.
 On the other hand, using the facts that $x_3\rightarrow y_2\rightarrow x_2$ (Claim 6(c)) and $x_2\rightarrow y_3$ we see that $Q:=x_1yx_3y_2x_2y_3\ldots x_1$ is a cycle of length $2a-4$. 
Notice that $x_1\rightarrow y_1\rightarrow x_2$ is a $Q$-bypass, whose the gap with respect to $Q$ is  equal to 4. This contradicts the minimality of  the gap $|C'|+1$ or Claim 6(a) since $d(y_1)\geq 2a-2$. Therefore, $d(x_2)\leq 2a-3$.  Together with condition $B_0$ this implies that $d(x_3)\geq 2a-2$ since $\{x_2,x_3\}\rightarrow y_2$. \fbox \\\\

\textbf{(b).} Suppose, on the contrary, that $|C'|=3$ and $d(y_1)\leq 2a-3$. By condition $B_0$,  $y_1 x_3\notin A(D)$ since 
$d(y)\leq 2a-3$ (Claim 6(a)) and $y\rightarrow x_3$. 
From (20) we have that  $y_0x_3\notin A(D)$ since  $yx_3\in A(D)$. 
Thus, the arcs  $y_1 x_3$ and  $y_0 x_3$ are not in $A(D)$.
This together with $d(x_3)\geq 2a-2$ (Claim 8(a)) imply that $x_3\rightarrow y_1$ and $y_3\rightarrow x_3$. 
By (ii),  $y_2\rightarrow x_4$ and hence,
 $Q:=x_1yx_3y_1x_2y_2x_4\ldots x_1$ is a cycle of length $2a-4$, which does not contain the vertices $x_0,y_0, x$ and $y_3$. The path $x_3y_3x_4$ is a $Q$-bypass whose the gap with respect to $Q$ is equal to 4. 
Therefore, by the minimality of the gap $|C'|+1$ and Claim 6(a), $d(y_3)\leq 2a-3$ but this is a contradiction since 
$\{y,y_3\}\rightarrow x_3$ and  
$d(y)\leq 2a-3$ (Claim 6(a)).  \fbox \\\\

\textbf{(c).} Suppose that the claim is not true. Without loss of generality, assume that for some $k\in [3,m]$, $x\rightarrow y_k$ and
$x_k\rightarrow y_1$. Then, by (ii) we have that $y_2\rightarrow \{x,x_0\}$ and hence, the cycle 
$x_1yx_3\ldots x_ky_1x_2y_2xy_k$ $\ldots  y_mx_1$ is a cycle of length $2a-2$, a contradiction. \fbox \\\\ 

Now we are ready to complete the proof of Theorem 1.7.
 Combining (19) and Claim 8(a), we obtain
$$
max\{d(x),d(x_0),d(y),d(y_0),d(x_2)\}\leq 2a-3. \eqno (21)
$$
This and condition $B_0$ imply 
$$
d^-(v,\{x,x_0,x_2\})\leq 1 \quad \hbox{for all} \quad v\in Y. 
$$
Therefore, since $d(y_1)\geq 2a-2$ (Claim 8(b)),  the vertex $y_1$ and every vertex of $X\setminus \{x,x_0,x_2\}$ form a 2-cycle, i.e.,
 $$y_1\leftrightarrow x_i \quad \hbox{for all} \quad i\in [1,m]\setminus \{2\}. \eqno (22)
$$
 Now using Claim 8(c), we obtain
$A(\{x,x_0\}\rightarrow \{y_2,y_3,\ldots , y_m\})=\emptyset$. 
From  $y_2 \mapsto \{x,x_0\}$ (Claim 6(d)) and the minimality of the gap $|C'|+1$ it follows that $d^-(y,\{x,x_0\})=0$. The last two equalities imply that 
$$
A(\{x,x_0\}\rightarrow \{y,y_2,y_3,\ldots , y_m\})=\emptyset. \eqno (23)
$$
By (20), in particular, we have
$$
max\{d^-(y_1,\{x,x_0\}),d^-(y_0,\{x,x_0\})\}\leq 1. \eqno (24)
$$
Since $D$ is strong, from (23) and (24) it follows  that $x\rightarrow y_0$ or $x_0\rightarrow y_0$. Again using  (24), we obtain that if 
 $x\rightarrow y_0$, then  $x_0 y_0\notin A(D)$ and  $x_0\rightarrow y_1$; if  $x_0\rightarrow y_0$, then  $x y_0\notin A(D)$ and  $x\rightarrow y_1$.

Because of the symmetry between the vertices $x$ and $x_0$, we can  assume that $x_0\rightarrow y_0$, $x\rightarrow y_1$ and $x_0 y_1\notin A(D)$, $x y_0\notin A(D)$. 
It is not difficult to show that $x_2$ and every vertex $y_i$ with $i\in [3,m]$ are not adjacent. Indeed, if 
$y_i\rightarrow x_2$, then, by (22), $y_1\rightarrow x_{i+1}$ and hence, $x_1yx_3y_3\ldots y_ix_2y_2xy_1x_{i+1}\ldots x_1$ is a cycle of length $2a-2$; 
 if 
$x_2\rightarrow y_i$, then, by (ii), $x_i\rightarrow y_2$ and hence, $x_1yx_3y_3\ldots x_iy_2xy_1x_2y_i\ldots y_mx_1$ is a cycle of length $2a-2$. Thus, in both cases we have a contradiction.

Therefore, $a(y_i,x_2)=0$ for all $i\in [3,m]$. This and (20) imply that $d(y_i)\leq 2a-3$ for all  $i\in [3,m]$.
From $y\rightarrow x_3$, $d(y)\leq 2a-3$ (Claim 6(a)) and condition $B_0$ it follows that 
$d^-(x_3,\{y_0,y_3,y_4,\ldots y_m\})=0$. 
Now frome  $d(x_3)\geq 2a-2$ (Claim 8(a)), we have that $m=3$, i.e., $a=5$.
 Since $x\rightarrow y_1$, $x_0\rightarrow y_0$ and  (21), it follows that $d^+(x_2,\{y_0,y_1\})=0$. Because of $d(y_1)\geq 2a-2$ (Claim 8(b)) and $d^-(y_1,\{x_0,x_2\})=0$ we have $y_1\rightarrow x_0$, $x_3\leftrightarrow y_1$ and $y_1\rightarrow x_1$.
 From this it follows that $y_0x_2\notin A(D)$ (for otherwise, $y_1x_0y_0x_2y_2\ldots x_1y_1$ is a cycle of length $2a-2$, a contradiction). Therefore, $a(x_2,y_0)=0$. If $y_0\rightarrow x_1$, then the cycle $x_1yx_3y_1x_2y_2x_0y_0x_1$ is a cycle of length $2a-2=8$, a contradiction. Therefore, $y_0 x_1\notin A(D)$. So, we have $d^+(y_0,\{x_1,x_2,x_3\})=0$. Then $d^+(y_0,\{x,x_0\})\geq 1$ since $D$ is strong. It is easy to see that $y_0x\notin A(D)$ (for otherwise, $x_1yx_3y_2x_0y_0xy_1x_1$ is a cycle of length 8, a contradiction. Therefore, $d^+(y_0,\{x_1,x_2,x_3, x\})= 0$.
 On the other hand, from (23) and $x_0y_1\notin A(D)$ we have  $N^+(x_0)=\{y_0\}$. Now it is not difficult to see that there is no path from $x_0$ to any vertex of $V(C)$ since $N^+(y_0)=\{x_0\}$, which contradicts that $D$ is strong. So, the discussion of the case $|C'|\geq 3$ is completed. Theorem 1.7 is proved. \\

\section {Concluding remarks}

In view of Theorem 1.3 it is natural to set the following problem.

\textbf{Problem.} {\it Characterize those strongly connected  balanced bipartite digraphs of order $2a\geq 6$ in which 
$d(x)+d(y)\geq 3a$ for every pair of vertices $x$, $y$ with a common in-neighbour or a common out-neighbour but are not even pancyclic.}

\end{document}